\documentclass{agtart_a}
\pdfoutput=1
\usepackage{amsmath}
\usepackage{xypic}

\title{Generalised Swan modules and the D(2) problem}
\author{Tim Edwards}
\givenname{Tim}
\surname{Edwards}

\address{Department of Mathematics\\
University College London\\\newline
Gower St\\
London WC1E 6BT\\UK}
\email{timeds@math.ucl.ac.uk}



\received{16 September 2005}
\revised{13 December 2005}
\accepted{9 January 2006}
\proposed{}
\seconded{}
\published{24 February 2006}
\publishedonline{24 February 2006}

\volumenumber{6}
\issuenumber{}
\publicationyear{2006}
\papernumber{2}
\startpage{71}
\endpage{89}

\doi{}
\MR{}
\Zbl{}

\makeatletter
\def\cnewtheorem#1[#2]#3{\newtheorem{#1}{#3}[section]
\expandafter\let\csname c@#1\endcsname\c@subsection}

\makeatother  

\def\overM{\wwbar{M}}
\def\overE{\wwbar{E}}
\def\tildeX{\wtilde{X}}
\def\SNF{{\rm SNF}}
\def\Diag{{\rm Diag}}

\newcommand{\lr}{\longrightarrow}
\newcommand{\equas}{\begin{eqnarray*}}
\newcommand{\equaf}{\end{eqnarray*}}
\newcommand{\End}{{\rm End}}
\newcommand{\Ext}{{\rm Ext}}

\newcommand{\Ker}{{\rm Ker}}

\newcommand{\zg}{\mathbb{Z}[G]}

\newcommand{\zG}{\mathbb{Z}[\Gamma]}
\newcommand{\zc}{\mathbb{Z}[C_{\infty}]}

\newcommand{\bem}{\left( \begin{array}{ccccccccc}}
\newcommand{\enm}{\end{array} \right)}
\newcommand{\bev}{\left( \begin{array}{c}}
\newcommand{\env}{\end{array} \right)}
\newcommand{\Mm}{\mathcal{M}_k(R_n)}

\cnewtheorem{conj}[subsection]{Conjecture}
\cnewtheorem{thm}[subsection]{Theorem}
\cnewtheorem{prop}[subsection]{Proposition}
\cnewtheorem{lem}[subsection]{Lemma}
\cnewtheorem{cor}[subsection]{Corollary}
\cnewtheorem{sublem}[subsubsection]{Lemma}
\cnewtheorem{subthm}[subsubsection]{Theorem}
\cnewtheorem{subprop}[subsubsection]{Proposition}
\cnewtheorem{subcor}[subsubsection]{Corollary}
\theoremstyle{definition}
\cnewtheorem{defi}[subsection]{Definition}
\cnewtheorem{cons}[subsection]{Construction}
\cnewtheorem{remark}[subsection]{Remark} 

\begin{document}

\begin{asciiabstract}
We give a detailed proof that, for any natural number n, each
algebraic two complex over C_n \times C_\infty is realised up to
congruence by a geometric complex arising from a presentation for the
group.
\end{asciiabstract}

\begin{abstract}
We give a detailed proof that, for any natural number $n$, each
algebraic two complex over $C_n{\times}C_{\infty}$ is realised up to
congruence by a geometric complex arising from a presentation for the
group. \end{abstract}

\maketitle
\section{Introduction}

It is well known that Whitehead's theorem allows the study of homotopy
types of two dimensional CW complexes to be phrased in terms of chain
homotopy types of algebraic complexes, arising as the cellular chains
on the universal cover. It is natural to ask whether the category of
algebraic complexes fully represents the category of CW complexes, in
particular whether every algebraic complex is realised
geometrically. The case of two dimensional complexes is of special
interest, partly due to the relationship between such complexes and
group presentations and partly since, as was recently proved, it
relates to the question as to when cohomology is a suitable indicator
of dimension.

By an {\it algebraic 2--complex} over a group $\Gamma$ we mean any
exact sequence of (right) $\zG$ modules of the form: $$ \xymatrix{
0 \ar[r] & M \ar[r] & F_2 \ar[r] & F_1 \ar[r] & F_0 \ar[r] & \Z \ar[r]
& 0, } $$ where each $F_i$ is finitely generated free and $\Z$
denotes the trivial $\zG$ module. The Realisation Problem asks if all
such complexes are chain homotopy equivalent to a complex arising from
a two dimensional CW complex, which may be assumed to be the Cayley
complex of some presentation for $\Gamma$. If so, we say that the {\sl
realisation property} holds for $\Gamma$.

 The D(2) Problem, as originally formulated by Wall in \cite{fc}, asks
if any three dimensional CW complex is necessarily homotopic to a
complex of dimension two provided that, ranging over all possible
coefficient systems, the complex has zero homology and cohomology in
dimensions higher than two. The problem is parameterised by the
fundamental group in the sense that, since homotopy equivalence
induces an isomorphism on fundamental groups, one may prove or
disprove the D(2) problem for CW complexes with a specified
fundamental group. Accordingly, we say that \textit{the D(2) property}
holds for $\Gamma$ if every three dimensional complex with fundamental
group isomorphic to $\Gamma$, satisfying the hypothesis of the D(2)
problem, is homotopy equivalent to a two dimensional complex.  In
\cite{SMD2}, Johnson relates the D(2) problem to the realization
problem as follows:

\begin{thm}[F\,E\,A Johnson] \label{ba} Let $\Gamma$ be a finitely presented group such that there is an algebraic 2--complex: 
$$0 \lr M \lr F_2 {\lr} F_1 {\lr} F_0 {\lr} \Z \lr 0$$ with $F_i$, $M$
finitely generated over $\zG$ and each $F_i$ free. Then the D(2)
property holds for $\Gamma$ if and only if the realization property
holds for $\Gamma$.
\end{thm}
We prove:

\medskip
{\bf Theorem A}\qua \textsl{The D(2) Property holds for $C_n \times
C_{\infty}$.}

\medskip
In the case of $\Gamma$ abelian and of rank one, the Bass--Murthy
paper \cite{Bass} shows that all stably free $\zG$ modules are
necessarily free. This result encouraged us to investigate the
realization problem in this case and we remark that the Bass--Murthy
result is essential to our work, although its role is subtle. Indeed,
in all confirmed cases of the D(2) property, the proof is dependent on
a proof that stably free modules are free. At present, the class of
groups to which this applies is limited to finite abelian groups
(Latiolais and Browning \cite{lat,browning}), free groups (Johnson
\cite{SMD2}), the dihedral groups of order $4n +2$ (Johnson
\cite{Min}) and the dihedral group of order 8 (Mannan
\cite{wajid}). In the appendix to this paper we give a proof that the
D(2) property holds for the group $C_{\infty} \times C_{\infty}$.

\section{Preliminaries}

Suppose that $\mathcal{G}= \langle x_1, \ldots x_g \mid W_1, \ldots ,W_r
\rangle$ is a presentation for $\Gamma$, then the corresponding Cayley
complex $X_{\mathcal{G}}$ is a 2--dimensional CW complex with
fundamental group $\Gamma$. The chain complex of the universal cover
$\tildeX_{\mathcal{G}}$ gives rise to a complex of $\zG$ modules thus:
$$C_*(\tildeX_{\mathcal{G}}) = \left( 0 \lr \pi_2( X_{\mathcal{G}})
\lr \zG^r \stackrel{\partial_2}{\lr} \zG^g \stackrel{\partial_1}{\lr}
\zG \stackrel{\varepsilon}{\lr} \Z \rightarrow 0 \right)$$ where:

\begin{itemize} 

\item[(i)] We have identified
Ker$(\partial_2)=H_2(\tildeX_{\mathcal{G}})$ with
$\pi_2(X_{\mathcal{G}})$ via the Hurewitz isomorphism.

\item[(ii)] Since the universal cover is simply connected, and by (i)
above, the complex is an exact sequence, so
$C_*(\tildeX_{\mathcal{G}})$ is an algebraic 2--complex.

\end{itemize} 

Note that $\partial_2$ and $\partial_1$ are completely determined by
the relations and generators of $\mathcal{G}$ and may be given
explicitly, following Fox \cite{Fox}.  In what follows, we shall say
that an algebraic 2--complex $E$ is {\sl geometrically realized} if is
there exists a presentation $\mathcal{G}$ for $\Gamma$ and a chain
homotopy equivalence $E \simeq C_*(\tildeX_{\mathcal{G}}) $. The
realization property, as defined in the introduction, holds for
$\Gamma$ if and only if every algebraic 2--complex over $\Gamma$ is
geometrically realized.

The first reduction is due to Schanuel \cite{Swanper} and shows that,
in the context of \fullref{ba}, $M$ is determined up to stability by
appearing in such a sequence -- where stable equivalence of modules
$M$ and $N$ is taken to be the existence of natural numbers $a$,$b$
and an isomorphism $M \oplus \zG^a \cong N \oplus \zG^b$. The first
task is to determine the class of modules which are stably isomorphic
to M, and has historically presented the more difficult objective. The
second determines each congruence class of algebraic two complexes
with a given $M$, and asks which of these are chain homotopy
equivalent to a complex arising from a group presentation. Two key
reductions are then possible in some cases:

\begin{prop}[F\,E\,A Johnson] \label{red1} 
Suppose that we are given $\Gamma$ satisfying the hypothesis of
\fullref{ba} and such that $\Ext^3(\Z,\zG)=0$. Then:

\begin{itemize} 
\item[(i)] 
If each algebraic 2--complex of the form $$ 0 \lr M \lr \zG^{a_1} \lr
\zG^{a_2} \lr \zG^{a_3} \lr \Z \lr 0$$ is geometrically realized, then
each algebraic 2--complex of the form $$ 0 \lr M \oplus \zG^{b_0} \lr
\zG^{b_1} \lr \zG^{b_2} \lr \zG^{b_3} \lr \Z \lr 0 $$ is geometrically
realized.
\item[(ii)] There exists a ring structure on $\textrm{Ext}^3(\Z, M)$
under which congruence classes of algebraic 2--complexes are
units. \end{itemize} \end{prop}

We present a concise version of the proof, following Humphreys
\cite{jodie}. Suppose that a presentation $\mathcal{H}$ for $\Gamma$
realizes the extension
$$C_{*}(\wtilde X_{\mathcal{H}}) = 0 \rightarrow M
\stackrel{i}{\rightarrow} \zG^{a_1} \stackrel{\delta_2}{\rightarrow}
\zG^{a_2} \stackrel{\delta_1}{\rightarrow} \zG^{a_3} \rightarrow \Z
\rightarrow 0,$$ then the addition of $n$ trivial relations to
$\mathcal{H}$ may be seen to realize the extension
$$f(C_{*}(\wtilde X_{\mathcal{H}}) ) = 0 \rightarrow M \oplus \zG^n
\stackrel{[i,Id]}{\lr} \zG^{a_1+n} \stackrel{[\delta_2,0]}{\lr}
\zG^{a_2} \stackrel{\delta_1}{\rightarrow} \zG^{a_3} \rightarrow \Z
\rightarrow 0.$$ If $E$ is an arbitrary extension in
$\textrm{Ext}^3(\Z,M)$, then we may define $f(E)$ to be the extension
in $\textrm{Ext}^3(\Z, M \oplus \zG^n)$ obtained similarly; $f$ is an
isomorphism $\textrm{Ext}^3(\Z, M) \cong \textrm{Ext}^3(\Z, M \oplus
\zG^n)$ and is a bijection on extensions which are chain homotopy
equivalent to algebraic complexes. This completes the proof of (i).

Define the relation $\sim$ on the endomorphism ring $\End (M)$ by $f
\sim h$ if and only if $f-h$ factors through $i\co M \rightarrow
\zG^{a_1}$. Then $\sim$ is an equivalence relation. By standard
homological algebra, there exists an additive group isomorphism
$$\Ext^3(\Z,M) \cong \End (M) / \sim \quad $$ we show that the ring
structure is preserved in the quotient, i.e.\ $ga \sim hb$ for any $a
\sim b$ and $g \sim h$. For any such $a,b,g,h \in \End(M)$ it is
trivial that $ga \sim gb$ (draw the diagram!). Since $(g-h) \sim 0$
there exists some factorisation $(g-h) = \varphi_1 \circ i$. Consider
the homomorphism $i \circ b\co M \rightarrow \zG^{a_1}$; this
represents a congruence class of $\Ext ^3(\Z,\zG^{a_1})$ and since
Ext$^3(\Z,\zG) = 0$, by additivity $i \circ b$ factors through
$\zG^{a_1}$ as some $\varphi_2 \circ i$ and $(g-h) \circ b = \varphi_1
\circ \varphi_2 \circ i$. Thus the product structure is well defined
on the quotient.

Suppose that some $f \in \End(M)$ represents the congruence class of
an algebraic 2--complex $E'$. Then $E'$ may be extended to a
projective resolution for $\Z$ which one may use to calculate the
elements of $\Ext^3(\Z ,M)$ as a quotient of End$(M)$. Let $h$
represent the congruence class of the given complex $E$ under these
terms, then $hf \sim Id$ and hence $h$ is a left inverse for
$f$. Since $\Ext^3(\Z ,M)$ is finitely generated as an abelian group,
$h$ is also a right inverse. The result follows.

 \section{Generalised Swan modules over $C_n \times C_{\infty}$}

Throughout, we shall assume that $G=C_n$ is the cyclic group on $n$
elements and $\Gamma$ is the product $G \times C_{\infty}$. Fix $R$ to
represent the ring and trivial $\zG$ module $\zc$, so that
$\zG=R[G]$. In this case, a standard representative of the stable
class of the second homotopy module may be constructed from the
presentation $$\mathcal{G} = \langle x,t | x^n=1,tx=xt \rangle.$$
Writing $N$ for $\sum _{g \in C_n} g$, $\pi_2( X_{\mathcal{G}})$ may
be identified with the submodule of $\zG^2$ generated by the elements
$$ \bev 0 \\ x-1 \env \quad {\rm and} \quad \bev N \\ t-1 \env.$$ with
corresponding chain complex \equas C_*(\tildeX_{\mathcal{G}}) = \left(
0 \lr \pi_2( X_{\mathcal{G}}) \lr \zG^2 \stackrel{\partial_2}{\lr}
\zG^2 \stackrel{\partial_1}{\lr} \zG \stackrel{\varepsilon}{\lr} \Z
\rightarrow 0 \right) \\ \partial_2 = \bem 1-t & N \\ x-1 & 0 \enm
\qquad \partial_1 = (x-1 \quad t-1). \qquad \qquad \equaf Since
$\pi_2( X_{\mathcal{G}})$ is finitely generated, $\Gamma$ satisfies
the hypothesis of Theorem 1.1. To see that $\Gamma$ satisfies the
conditions of \fullref{red1}, note first that the augmentation ideal
$I=\Ker(\varepsilon\co \zG \rightarrow \Z)$ is free as an $R$ module,
generated by the elements $(t-1)$ and $\{(x^i-1)\co1\leq i \leq
n-1\}$. The class of such modules is {\sl tame} in the sense of
Johnson \cite[Chapter 4]{SMD2}, so that projective modules are
relatively injective and $\Ext ^2(I, \zG)=0$. The result follows by
dimension shifting.

It shall prove easier to calculate the stable class of the dual module
$\pi_2( X_{\mathcal{G}})^*$, and we justify doing so by
\cite[Proposition 28.1]{SMD2}, which shows that cancellation holds for
$\pi_2(X_{\mathcal{G}})^*$ if and only if it holds for $\pi_2(
X_{\mathcal{G}})$.

Directly, one may verify that $\pi_2( X_{\mathcal{G}})^*$ may be
identified with the submodule of $\zG^2$ generated by the the
elements:
$$ \bev 1-t \\ x-1 \env \quad \bev N \\ 0 \env $$ $R$ imbeds in $\pi_2(
X_{\mathcal{G}})^*$ by identification with the module generated by the
second element above, resulting in a short exact sequence of $\zG$
modules: $$ 0 \rightarrow R \rightarrow \pi_2( X_{\mathcal{G}})^*
\rightarrow S \rightarrow 0$$ where $S = \zG / \langle N
\rangle$. Embed $R$ into $\zG$ by identifying $R$ with the submodule
generated by $N$, so that for each $k \in \mathbb{N}$ there is an
exact sequence: $$ \xymatrix{ 0 \ar[r] &R^k \ar[r] & \pi_2(
X_{\mathcal{G}})^* \oplus \zG^{k-1} \ar[r] & S^k \ar[r] & 0}$$
giving a na\"\i ve model for modules potentially stably isomorphic to
$\pi_2( X_{\mathcal{G}})^*$ as a subclass of the central modules
occurring in extensions of $R^k$ by $S^k$. Through a detailed study of
such extensions we shall show that this model is appropriate and prove
the required cancellation result.

\begin{prop} 
Let $\Z_n$ denote the integers modulo $n$, and $\mathcal {M}_k (R_n)$
the ring of $k$ by $k$ matrices over the ring $\Z_n[C_{\infty}$. Then
$\Ext ^1_{\zG}(S^k,R^k) \cong \Mm.$ \end{prop}

\begin{proof} The truth of this statement may be deduced from the `change of rings formula' given in Cartan and Eilenberg \cite{CEHom} and from the cohomology of $C_n$. We give two explicit proofs, with the justification that the given constructions of congruence classes will prove useful in the sequel.

A projective resolution of $S^k$ is given by: 
$$ 
\xymatrix{ \ldots
\ar[r]^{\hbox{\llap{$\scriptstyle(x-1)I_k$\hglue-8pt}}} & \zG^k \ar[r]^{(N)I_k} & \zG^k \ar[r]^{\pi} & S^k
\ar[r] & 0 } 
$$ 
where $S^k$ is identified with $\zG^k / (N)\zG^k$
and $I_k$ denotes the identity matrix. Note that the kernel of $\pi$
may be identified with $R^k$. Through the standard homological
classification of $\Ext $, Mac Lane \cite{Mac}, it is enough to
determine the additive group of homomorphisms $f\co R^k \rightarrow
R^k$, modulo those which factor through the inclusion of $R^k$ in
$\zG^k$. Since $R$ may be considered as a ring, and $f$ an $R$ module
homomorphism, we may represent $f$ as multiplication by some $k \times
k$ matrix $F$ with entries in $R$. Since the image of $R^k$ in $\zG^k$
contains only elements divisible by $N$, if $f$ factors through
$\zG^k$, then all of the entries in $F$ must be divisible by
$n$. Conversely if $F$ is divisible by $n$, then $f$ clearly factors
through $\zG^k$, i.e.\ $F$ represents a unique congruence class modulo
$(n)$. For each matrix $F$ in $\mathcal{M}_k(R_n)$, we may pick a lift
for $F$ in $\mathcal{M}_k(R)$ and define $E(F)$ to be the extension
constructed through the pushout: $$ \xymatrix{ 0 \ar[r] & R^k
\ar[r] \ar[d]^{F} & \zG^k \ar[r]^{\pi} \ar[d] & S^k \ar[r] \ar@{=}[d]&
0 \\ 0 \ar[r] & R^k \ar[r] & M(F) \ar[r] & S^k \ar[r] & 0 } $$
where the bottom row is $E(F)$ and $M(F)$, the central module, is
determined by $F$. The correspondence between $\Ext ^1_{\zG}(S^k,R^k)$
and $\Mm$ clearly commutes with the addition operations, and our first
characterisation is complete.

We perform a second calculation following the methods of \cite{SMD2}
regarding extensions within tame classes. Since $R,S$ and $\zG$ are
free as $R$ modules we may take a (relatively) injective coresolution
of $R^k$ to be: $$ \xymatrix{0 \ar[r] & R^k \ar[r] & {\,\zG^k\,}
\ar[r]^{(x-1)I_k} & {\,\zG^k\,} \ar[r]^{(N)I_k} &\ldots } $$ We wish to
determine the group of homomorphisms $f\co S^k \rightarrow S^k$ modulo
those which factor through $\pi$ as $f=\pi \circ h$ for some $h\co S^k
\rightarrow \zG ^k$. Since $S=\zG / (N)$ may be considered as a
quotient ring of $\zG$, each homomorphism $f\co S^k \rightarrow S^k$
is a $S$ module homomorphism and may be represented as multiplication
by some $k \times k$ matrix $F$ with entries in $S$. By observation of
the coresolution above, $f$ factors through $\pi$ if and only if
$(x-1)$ divides each entry of $F$. Note that $\hat \varepsilon$
induces a surjective augmentation homomorphism $\hat \varepsilon _S\co
S \rightarrow R_n$, which we extend to $\hat \varepsilon_{S}
\co\mathcal{M}_k(S) \rightarrow \Mm$. Then the kernel of $\hat
\varepsilon_S$ is precisely the set of matrices with $(x-1)$ a divisor
of each entry. This establishes a one-to-one correspondence between
$\Ext ^1_{\zG}(S^k,R^k)$ and $\Mm$ which, again, clearly commutes with
addition.  For each matrix $F$ in $\mathcal{M}_k(R_n)$ we may pick a
lift of $F$ in $\mathcal{M}_k(S)$ and define $\overE(F)$ to be the
bottom row in the pullback: $$ \xymatrix{ 0 \ar[r] & R^k \ar[r] &
\zG^k \ar[r]^{\pi} & S^k \ar[r] & 0 \\ 0 \ar[r] & R^k \ar@{=}[u]
\ar[r] & \overM(F) \ar[r] \ar[u] & S^k \ar[r] \ar[u]^{F} & 0 } $$
so that each extension of $R^k$ by $S^k$ is congruent to $\overE(F)$
for some unique $F \in \mathcal{M}_k (R_N)$.
\end{proof}
To give concrete examples of the congruence classes and central
modules, given any $A \in \Mm$, we may pick a lift of $A$ in
$\mathcal{M}_k(\zG)$ and define $M(A)$ to be the submodule of
$\zG^{2k}$ generated by the columns of the matrix: \equas A'= \bem A &
(N)I_k \\ (x-1)I_k & 0 \enm \equaf Then $R^k$ embeds in $M(A)$ by
identifying the $i^{th}$ basis element of $R^k$ with the $i+k^{th}$
column of $A'$, and the image of $R^k$ in $M(A)$ is the kernel of the
surjective map $\pi_A\co M(A) \rightarrow S^k$ given by sending the
$i^{th}$ column of $A'$ (for $1 \leq i \leq k)$ to the $i^{th}$
generator of $S^k$. We represent the resulting exact sequence as:
$$ \xymatrix{E'(A) = & 0 \ar[r] & R^k \ar[r]^{i_A} & M(A)
\ar[r]^{\pi_A} & S^k \ar[r] & 0 } $$ and $E'(A) \equiv E(A) \equiv
\overE (A)$.

We wish to characterise the isomorphism classes of the modules
$M(A)$. Since the action of $N$ is zero on $S$ and multiplication by
$n$ on $R$, and since $S$ is torsion free:
 $${\rm Hom}_{\zG}(R^k,S^k)=0,$$ so that any isomorphism $\phi\co M(A)
\cong M(B)$ induces an isomorphism of extensions $E(A) \cong E(B)$,
where the isomorphisms on each end may be represented as matrices $C$
and $D$. We distinguish the matrices in $\mathcal{M}_k(R_n)$ which are
images of such isomorphisms: \begin{itemize} \item $GR(k)$ denotes the
image of $GL_k(R)$ in $GL_k(R_n)$ under the natural map. \item $GS(k)$
denotes the image of $GL_k(S)$ in $GL_k(R_n)$ under
$\varepsilon_S$. \end{itemize}

Explicitly, any isomorphism $\phi\co M(A) \cong M(B)$ results in a
commutative diagram: $$ \xymatrix{ & 0 \ar[r] & R^k \ar[r]
\ar[d]^{C} & M(A) \ar[r] \ar[d]^{\phi} & S^k \ar[r] \ar[d]^D & 0 \\ &
0 \ar[r] & R^k \ar[r] & M(B) \ar[r] & S^k \ar[r] & 0 } $$ and it
is easily shown that the image of $C^{-1}AD$ in $\Mm$ must be equal to
$B$, so that $M(B) \cong M(C^{-1}AD)$. A diagram similar to the above
gives, through the five lemma:

\begin{thm}\label{iso} For arbitrary matrices $A,B \in \Mm$, $M(A) \cong M(B)$ if and only if there exists $C \in GR(k), D \in GS(k)$ such that $CAD = B$.
\end{thm}

Having characterised the isomorphism classes of the modules $M(A)$, we
wish to characterise their stable isomorphism classes. Clearly $$M(A)
\oplus \zG^m \cong M(A \oplus I_m),$$ where for $B \in
\mathcal{M}_m(R_n)$, the matrix $A \oplus B$ is taken to be the $(k+m)
\times (k+m)$ matrix: \equas \bem A & 0 \\ 0 & B \enm \equaf

\begin{cor} $M(A)$ is stably equivalent to $M(B)$ if and only if there exists $C \in GR(k+m), D \in GS(k+m)$ such that $C(A \oplus I_m)D = (B \oplus I_m)$. \end{cor}
As is standard, we write $E_k(R_n)$ for the group of $k \times k$
matrices over $R_n$ which are products of elementary matrices. Note
that $E_k(R_n) \subseteq GS(k)$ and $E_k(R_n) \subseteq GR(k)$. The
reader may wish to recall some of the basic results in algebraic
$K$--theory, available for example in Milnor \cite{Milnor}, from which
we take our notation.
\begin{prop}\label{det} The determinant homomorphism ${\rm det}\co GL_k(R_n) / E_k(R_n) \rightarrow R_n^+$ is an isomorphism for all $k \geq 1$. Equivalently, every invertible matrix in $\Mm$ with determinant one is a product of elementary matrices. \end{prop}

\begin{proof} Let the prime decomposition of $n$ be $\Pi_{i=1}^{s}p_i^{e_i}$. Then $R_n$ is isomorphic to the product $\Pi_{i=1}^{s}R_{p_i^{e_i}}$ and $\Mm$, $R_n^+$ decompose as products similarly, as does the determinant homomorphism. A fortiori, it is enough to prove the proposition in the case where $n = p^e$ is a power of a prime. 

Suppose that $\det (E) = 1$ where $E \in \Mm$, and let $E_p$ denote
the equivalence class of $E$ in $\mathcal{M}_k(R_p)$. We shall show
that we may reduce $E$ by elementary row and column operations to the
identity. Since $R_p=\mathbb{F}_p[t,t^{-1}]$ is a Euclidean domain
there are elementary matrices $E^1_p, \ldots E^j_p$ such that $E_p
\cdot E^1_p \cdot \ldots E^j_p= Id_p$. Thus we may pick elementary
matrices $E^1, \ldots E^j \in E_k(R_n)$ such that: \equas E \cdot E^1
\ldots E^j = \bem 1 + c_{1,1}p & c_{1,2}p & \ldots & c_{1,k}p \\
c_{1,2}p & 1+c_{2,2}p & \ldots & c_{2,k}p \\ \vdots & \vdots & \quad &
\vdots \\ c_{k,1}p & c_{k,2}p & \ldots & 1+c_{k,k}p \enm \qquad
c_{i,j} \in R_n\equaf By the binomial theorem each $(1+c_{i,i}p)$ is a
unit (take the $(p^e)^{th}$ power!) and thus the above matrix has
diagonal entries which are all units. Thus the above matrix, and hence
$E$, may be reduced by the action of column operations to a matrix of
the form: \equas \bem 1 + d_{1,1}p & 0 & \ldots & 0 \\ 0 & 1+d_{2,2}p &
\ldots & 0 \\ \vdots & \quad & \quad & \vdots \\ 0 & 0 & \ldots &
1+d_{k,k}p \enm \qquad d_{i,i} \in R_n \equaf By Whitehead's Lemma,
$\Diag(1,\ldots,1,u,u^{-1},1,\dots,1) \in E_k(R_n)$ and $E$ may be
further reduced to $$E'=\Diag(u,1,\ldots,1)$$ for some $u \in
R_n^+$. Since all elementary matrices have determinant one, ${\rm
det}(E') = {\rm det}(E)=1$, i.e.\ $u=1$ and $E'=I_k$. Hence $E$ is a
product of elementary matrices.
\end{proof}

Recall that for any Euclidean domain $\mathcal{R}$ and any $k \times
k$ matrix $A$ over $\mathcal{R}$, there are elementary matrices $E_1$
and $E_2$ over $\mathcal{R}$ such that \begin{itemize} \item $E_1\cdot
A\cdot E_2=\Diag(a_1,a_2,\ldots,a_k)$ is a diagonal matrix with
$|a_{i+1}|$ a divisor of $|a_{i}|$. \item $\Diag(a_1,a_2,\ldots,a_k)$
is unique up to multiplication by diagonal elementary matrices and is
sometimes called the {\sl Smith Normal Form} of $A$, denoted
$\SNF(A)$. \end{itemize} $E_1$ and $E_2$ may easily be constructed
through a process derived from the Gauss algorithm for Euclidean
domains.
\begin{thm}\label{cancdet} For each non-zero $\alpha \in R_n$ and all $B \in \Mm$, $M(\alpha)$ is stably equivalent to $M(B)$ and only if $M(B) \cong M(\alpha) \oplus \zG^{k-1}$. \end{thm}

\begin{proof} Suppose that $M(\alpha)$ is stably equivalent to $M(B)$ for some non-zero $\alpha \in R_n$ and some $B \in \Mm$. Then by $\zc$ rank considerations we may assume that there exists an $m \in \mathbb{N}$ and an isomorphism $M(B) \oplus \zG^m \cong M(\alpha) \oplus \zG^{k+m-1}$, so that by \fullref{iso} there exist matrices $C \in GR(k+m)$, $D \in GS(k+m)$ such that $$C(B \oplus I_m)D = (\alpha \oplus I_{m+k-1}).$$
Define $ B_{new} = (\det (C) \oplus I_{k-1})\cdot B \cdot (\det (D)
\oplus I_{k-1})$. Since $R$ and $S$ are commutative rings, the
determinant homomorphisms are well defined and \equas (\det (C) \oplus
I_{k-1}) \in GR(k) \\ (\det (D) \oplus I_{k-1}) \in GS(k).\equaf By
\fullref{iso} $M(B) \cong M(B_{new})$, so that we may assume that
$B=B_{new}$ and $\det(B) = \alpha$. We shall show that $B$ must be
reducible by the action of elementary matrices to $(\alpha \oplus
I_{k-1})$, which would imply the result. Again, this is essentially a
statement about matrices over $R_n$ and we may assume that $n=p^{e}$
is a power of a prime. As in the proof of the lemma, let $B_p$ denote
the congruence class of $B$ mod $p$. Since $R_p$ is a Euclidean
domain, we may reduce $B_p$ by row and column operations to a matrix
of the form $\SNF(B_p)=\Diag ( b_1, \ldots b_k)$ and moreover we we
may assume that $\SNF(B_p) = \Diag(\alpha_p,1, \ldots 1)$. Thus, over
$R_n$ with $n=p^e$, we may reduce $B$ by the action of elementary
matrices to a matrix of the form: \equas \bem \alpha + c_{1,1}p &
c_{1,2}p & \ldots & c_{1,k}p \\ c_{1,2}p & 1+c_{2,2}p & \ldots &
c_{2,k}p \\ \vdots & \quad & \quad & \vdots \\ c_{k,1}p & c_{k,2}p &
\ldots & 1+c_{k,k}p \enm \qquad c_{i,j} \in R_n\equaf Again, each
$1+c_{i,i}p$ is a unit and we may reduce $B$ to a matrix of the form
$\Diag(\alpha + dp, 1, \ldots ,1)$ for some $d \in R_n$. This
completes the proof, since we have insisted that ${\rm det
}(B)=\alpha$ and hence $B$ may be reduced by the action of elementary
matrices to $(\alpha \oplus I_{k-1})$.
\end{proof}

Theorem 3.5 shows that a limited form of cancellation holds within the
class of modules $M(A)$, of which $\pi_2( X_{\mathcal{G}})^* = M(t-1)$
is a member. In order to extend this class we shall employ a technique
for constructing projectives over fibre products due to Milnor
\cite{Milnor}. A {\sl Cartesian} square of rings is a commutative
diagram of rings and ring homomorphisms: $$ \xymatrix{ \mathcal{R}
\ar[r]^{i_1} \ar[d]^{i_2} & \mathcal{R}_1 \ar[d]^{j_1} \\
\mathcal{R}_2 \ar[r]^{j_2} & \mathcal{R}_0 } $$ such that
 
 \begin{itemize} \item[(i)] $\mathcal{R}$ is the fibre product of
 $\mathcal{R}_1$ and $\mathcal{R}_2$ over $\mathcal{R}_0$. This means
 that $\mathcal{R}$ may be identified with the set of pairs $(r_1,r_2)
 \in \mathcal{R}_1 \times \mathcal{R}_2$ such that
 $j_1(r_1)=j_2(r_2)$. \item[(ii)] At least one of $j_1$ or $j_2$ is
 surjective. \end{itemize}
 
If $P_i$ is a projective module over $\mathcal{R}_i$ then there is an
induced $\mathcal{R}_0$ module $j_{i \#} (P_i)$ given by $P_i
\otimes_{\mathcal{R}_i} \mathcal{R}_0$ and a canonical $R_i$--linear
map $j_{i *} \co P_i \rightarrow j_{i \#} (P_i)$ given by $$j_{i*} (p)
= 1 \otimes p.$$ Given projective modules $P_1$ over $\mathcal{R}_1$,
$P_2$ over $\mathcal{R}_2$ and an $\mathcal{R}_0$ isomorphism
$$h\co j_{1 \#} (P_1) \cong j_{2 \#} (P_2),$$ let $M(P_1,P_2,h)$
denote the subgroup of $P_1 \times P_1$ consisting of all pairs
$(p_1,p_2)$ such that $hj_{1*}(p_1) = j_{2*}(p_2)$. Then
$M(P_1,P_2,h)$ has a natural $\mathcal{R}$ module structure.

\begin{prop}[J Milnor]\label{Milproj}The module $M(P_1,P_2,h)$ is projective over $\mathcal{R}$, and every projective $\mathcal{R}$ module is isomorphic to some $M(P_1,P_2,h)$. Moreover, the modules $P_1$ and $P_2$ are naturally isomorphic to $i_{1 \#} M(P_1,P_2,h)$ and $i_{2 \#} M(P_1,P_2,h)$ respectively. \end{prop} 
We may now prove the following theorem.

\begin{thm}\label{canc} If $M$ is any module such that $M \oplus \zG^m \cong M(B)$ for some $B \in \mathcal{M} _{k+m}(R_n)$, then $M \cong M(A)$ for some $A \in \Mm$. \end{thm} 

We remark that it is here that we use the work of Bass--Murthy. In
particular we shall use the details in section 9 of \cite{Bass}.
\begin{proof} It is sufficient to show that there exists an exact sequence
 \equas 0 \rightarrow R^k \rightarrow M \rightarrow S^k \rightarrow 0
 \equaf where as before $S = \zG / (N)$. It is clear that, in the
 extension $E(B)$, $R^{k+m}$ is identified precisely with the
 submodule of $M(B)$ on which the action of $x$ is trivial.  Let $M_G$
 denote the submodule of $M$ on which $x$ acts trivially, so that
 there is an exact sequence:$$ 0 \rightarrow M_G \rightarrow M
 \rightarrow M / M_G \rightarrow 0.$$ Any isomorphism $M \oplus \zG^m
 \cong M(B)$ induces isomorphisms $M_G \oplus R^m \cong R^{k+m}$, and
 $M / M_G \oplus S^{m} \cong S^{k+m}$. Thus $M_G$ is stably free, and
 hence free, as an $R$ module. Since the action of $x$ on $M_G$ is
 trivial, we deduce that $M_G \cong R^k$ as a $\zG$ module. It remains
 to show that $M / M_G \cong S^k$. Note that $M / M_G$ has an $S$
 module structure which is then stably free, and if we can show that
 $M / M_G$ is free over $S$ we may deduce the result.

We claim that all stably free $S$ modules are free. Let $S_{\Z}$
denote $\zg / (N)\zg$, where as before $G = C_n$, so that there is a
natural identification of $S$ with $S_{\Z}[C_{\infty}$. By
Bass--Murthy \cite[9.1]{Bass} it is enough to prove that $S_{\Z}$ has
finitely many non-projective maximal ideals. We shall deduce it from
the fact that $\zg$ has finitely many non-projective maximal ideals,
(which \cite{Bass} claims is ``not difficult to verify", and we prove
for the sake of completeness in an appendix).

Let $i_1 \co \zg \rightarrow S_{\Z}$ be the natural map onto the
quotient, $i_2 \co \zg \rightarrow \Z$ be augmentation, $j_1\co S_{\Z}
\rightarrow \Z_n$ be defined as $j_1(\alpha + (N)) = i_2(\alpha) \mod
n$ and $j_2 \co \Z \rightarrow \Z_n$ be the natural map. Then the
following is a commutative diagram of rings and surjective ring
homomorphisms: $$ \xymatrix{ \zg \ar[r]^{i_1} \ar[d]^{i_2} &
S_{\Z} \ar[d]^{j_1} \\ \Z \ar[r]^{j_2} & \Z_{n} } $$ and $\zg$ may
be identified with the fibre product of $S_{\Z}$ and $\Z$ over
$\Z_n$. Let $I$ be a maximal ideal of $S_{\Z}$, and let $I'$ denote
the (necessarily projective) maximal ideal of $\Z$ given by $I' =
j_2^{-1}(j_1(I))$. Then $(I,I')$ is a maximal ideal of $\zg$ which is
isomorphic to the module $M(I,I',Id)$ constructed as
above. Furthermore, by \fullref{Milproj}, (I,I') is projective if and
only if $I$ is projective. Lastly, if $J$ is another maximal ideal of
$S_{\Z}$, then $(J,J') = (I,I')$ if and only if $J=I$ and hence there
is an injective map from the maximal ideals of $S_{\Z}$ to the maximal
ideals of $\zg$ which preserves projectivity.
\end{proof}

\begin{cor}\label{cancfinal} 
Any module stably isomorphic to $\pi_2( X_{\mathcal{G}})$ is
necessarily isomorphic to $\pi_2( X_{\mathcal{G}}) \oplus \zG^m$ for
some $m$, where $\pi_2( X_{\mathcal{G}})$ is the module defined in the
beginning of this section.
\end{cor}

\begin{proof} Observe, from our description of the generators of $\pi_2(
X_{\mathcal{G}})^*$, that there is an isomorphism $\pi_2(
X_{\mathcal{G}})^* \cong M(t-1)$. The required result is then a clear
consequence of \fullref{cancdet} and \fullref{canc}.
\end{proof}

\section{The D(2) problem for $C_n \times C_{\infty}$} 

We bring together the work of the previous sections:
\begin{prop} 
In order to prove the D(2) Problem for $C_n \times C_{\infty}$ it is
sufficient to realise geometrically all algebraic 2--complexes of the
form:
$$ 0 \lr \pi_2( X_{\mathcal{G}}) \lr \zG^a \lr \zG^b \lr \zG^c \lr \Z
\lr 0. $$ \end{prop}
\begin{proof} By \fullref{ba} we may consider the D(2) problem to be equivalent to the realization problem for $\Gamma$. The result then follows from \fullref{red1} and from \fullref{cancfinal}. \end{proof}

\begin{lem}\label{bd} There is an isomorphism {\rm Ext}$^3 _{\zG} (\Z,\pi_2( X_{\mathcal{G}})) \cong \Z_n$.\end{lem}
 
\begin{proof} One may calculate this result directly using the resolution: 
$$ \cdots \stackrel{\partial_2}{\longrightarrow}\zG^2
\stackrel{\partial_3}{\longrightarrow} \zG^2
\stackrel{\partial_2}{\longrightarrow} \zG^2
\stackrel{\partial_1}{\longrightarrow} \zG \lr \Z \lr 0$$ where
$\partial_1,\partial_2$ are as before and $$\partial_3 = \bem 0 & N \\
x-1 & t-1 \enm.$$ We simplify our proof by noticing that, immediately
from the generators given for $ \pi_2( X_{\mathcal{G}})$ in the start
of the previous section, there is a clear isomorphism $ \pi_2(
X_{\mathcal{G}})\cong I$, where $I$ is the augmentation ideal and is
generated as a submodule of $\zG$ by the elements $(x-1)$ and $(t-1)$.

If $f\co\zG^2 \lr I$ is a homomorphism, then we may represent $f$ as
multiplication by some matrix
$$\left( a \quad b \right) \qquad a,b \in I$$ $f$ is a cocycle if and
only if $a=\alpha_1(x-1)$, $b=(\alpha_1+\alpha_2 N)(t-1)$ for some
$\alpha_1 \in \zG, \alpha_2 \in R$; $\alpha_1$ is well defined modulo
$(N)$ and $\alpha_2$ is determined by $\alpha_1$.

$f$ is a coboundary if there exists some $a,b \in I$ such that
$$(a \quad b) \cdot \bem 0 & N \\ x-1 & t-1 \enm = \Big(
\alpha_1(x-1) \qquad (\alpha_1 + \alpha_2N)(t-1)\Big),$$ from which
we deduce that $f$ if a coboundary if and only if there exists a
choice of $\alpha_1 \in I$. The reader may verify that the map
$\varphi\co{\rm Ext}^3 _{\zG} (\Z,\pi_2( X_{\mathcal{G}})) \rightarrow
\Z_n$ given by:
$$\varphi \left( \alpha_1(x-1) \qquad (\alpha_1 +
\alpha_2N)(t-1)\right) \quad =\quad \varepsilon(\alpha_1) \mod n$$ is
a (ring) isomorphism.  \end{proof}
 
\begin{proof}[Proof of Theorem A] We have already shown that it is sufficient to realize all extensions of the form
$$ 0 \lr \pi_2( X_{\mathcal{G}}) \lr \zG^a \lr \zG^b \lr \zG^c \lr \Z
\lr 0 $$ and by \fullref{red1}, $\Ext^3(\Z,\pi_2( X_{\mathcal{G}}))$
has the structure of a ring under which algebraic 2--complexes are
necessarily units. For each unit $w \in \Z_n^+$ we shall realize the
class of $w$ up to congruence, and hence up to chain homotopy
equivalence.

An obvious change one may make to the standard presentation
$$\mathcal{G} = \left\langle\textrm{ } x,t \mid x^n=1,tx=xt\textrm{ } \right\rangle$$ is to replace the
generator $x$ for $C_n$ by the generator $y(v)=x^{v}$ where $1 \leq v
\leq n-1$ is a natural number coprime to $n$. Denote each such
presentation by
$$\mathcal{G}(v) = \left\langle\textrm{ } y(v),\textrm{ } t \textrm{ }
| \textrm{ } y(v)^n = 1 \textrm{ } , \textrm{ } y(v) t = t y(v)
\textrm{ } \right\rangle$$ where $y(v)=x^v$.

We remark that the Cayley complex of $\mathcal{G}(v)$ is homotopy
equivalent to the standard one, since as presentations they are
identical. Indeed, one may view this complex as that arising from
changing the isomorphism between fundamental group of the original
complex determined by $\mathcal{G}$ and the abstract group determined
by the presentation. The aforementioned homotopy equivalence does not
induce the identity on the fundamental groups of the spaces, and so
does not imply that the resulting algebraic complexes are chain
homotopic, although the distinction is meaningless from a geometric
viewpoint.

The corresponding chain complex of the universal cover of
$X_{\mathcal{G}(v)}$ is then \equas 0 \lr \Ker(\partial^v_2) \lr \zG^2
\stackrel{\partial^v_2}{\longrightarrow} \zG^2
\stackrel{\partial^v_1}{\longrightarrow} \zG
\stackrel{\varepsilon}{\longrightarrow} \Z \lr 0 \\ \textrm{where}
\quad \partial^v_2 = \bem1-t & N \\ x^v-1 & 0 \enm \qquad \textrm{and}
\qquad \partial^v_1 = (x^v-1 \quad t-1) \equaf so that
$\Ker(\partial^v_2) = \Ker(\partial_2)=\pi_2( X_{\mathcal{G}})$. Let
$w$ denote the inverse of $v$ mod $n$. Set $\tau$ to be the element $1
+ x^v + \ldots x^{v(w-1)}$, then $(1-x^v)\tau = (1-x)$ and the
following diagram commutes:
\begin{equation*} \xymatrix {0 \ar[r] & \pi_2( X_{\mathcal{G}}) \ar[r]^{j} \ar[d]^{f_3} &  \zG^2
\ar[r]^{\partial_2} \ar[d]^{f_2} & \zG^2 \ar[d]^{f_1}
\ar[r]^{\partial_1} & \zG \ar[d]^{Id} \ar[r]^{\varepsilon} & \Z
\ar[d]^{Id} \ar[r] & 0 \\ 0 \ar[r] & \pi_2( X_{\mathcal{G}})
\ar[r]^{j} & \zG^2 \ar[r]^{\partial^v_2} & \zG^2 \ar[r]^{\partial^v_1}
& \zG \ar[r]^{\varepsilon} & \Z \ar[r] & 0} \end{equation*} where each
$f_i$ may be represented as multiplication on the left by:
 $$f_3 = \tau\quad f_2 = {\bem {\tau} & 0 \\ 0 & \tau \enm}\quad f_1={\bem
 {\tau} & 0 \\ 0 & 1 \enm} $$ Finally, since $\varepsilon (\tau) = w
 $, $C_*(\tildeX_{\mathcal{G}(v)})$ realizes geometrically the
 congruence class of $w \in \Ext^3(\Z,\pi_2( X_{\mathcal{G}}))$. This
 completes the proof.
\end{proof}

 \section{Swan modules and weak cancellation} The machinery developed
to prove the cancellation result in section 3 may be applied to Swan
modules over a finite group, where $R=\zc$ is replaced by $\Z$.

Suppose that $G=\{g_i\}_{i=1}^{n}$ is a finite group of order $n$,
then the trivial $\zg$ module $\Z$ embeds into $\zg$ by identification
with the ideal generated by the element $N=\sum g_i$. There is a
resulting exact sequence:
$$\xymatrix{0 \ar[r] & \Z \ar[r] & \zg \ar[r] & \zg / (N) \ar[r] &
0.}$$ Similarly to before, we shall write $S$ for the $\zg$ module and
quotient ring $\zg / (N)$. Note that $S$ is the dual of the
augmentation ideal, and may be identified with the submodule of
$\zg^n$ generated by the element: $$\bem g_1 - 1 \\ \vdots \\ g_n - 1
\enm$$ For each $k \in \mathbb{N}$ and $A \in \mathcal{M}_k(\Z_n)$,
we may define $M(A)$ to be the submodule of $\zg^{k(n+1)}$ generated
by the columns of the matrix:
$$A' = \bem A'' & N \cdot I_k \\ (g_1 - 1) I_k & 0 \\ \vdots & \vdots
\\ (g_n - 1) I_k & 0 \enm$$ where $A''$ is a lift of $A$ under the
natural map $\zg \rightarrow \Z_n$. Then we may identify $\Z^k$ with
the submodule of $M(A)$ generated by the last $k$ columns of $A'$, and
the quotient under this embedding is clearly isomorphic to
$S^k$. Define $\varphi \co \zg^k \rightarrow M(A)$ by sending the
$i^{th}$ generator of $\zg$ to the $i^{th}$ column of $A'$, so that
the following commutes: $$ \xymatrix{ 0 \ar[r] & \Z^k \ar[r]
\ar[d]^{A} & \zg^k \ar[r] \ar[d]^{\varphi} & S^k \ar[r] \ar@{=}[d] & 0
\\ 0 \ar[r] & \Z^k \ar[r] & M(A) \ar[r] & S^k \ar[r] & 0 } $$ and
we may define $E(A)$ to be the congruence class of the bottom row. Let
$GR(k)$ denote the image of $GL_k(\Z)$ in $GL_k(\Z_n)$, and $GS(k)$
the image of $GL_k(S)$. Note that since $GR(k)$ contains matrices of
determinant $\pm 1$ we have $GR(k) \subseteq GS(k)$. The reader may
verify that the analysis of section 3 transfers immediately to such
extensions and proves the following:
\begin{prop}\label{simp} $\quad$
 \begin{itemize}
\item[(a)] Each extension of $\Ext^1(S^k,\Z^k)$ is congruent to $E(A)$
for some $A \in \mathcal{M}_k(\Z_n)$ and if $A_1$, $A_2$ are both
lifts of $A$, then $E(A_1) \equiv E(A_2)$.
\item[(b)] For arbitrary matrices $A,B \in \mathcal{M}_k(\Z_n)$, $M(A)
\cong M(B)$ if and only if there exists $C \in GR(k), D \in GS(k)$
such that $CAD = B$.
\item[(c)] For each $A \in \mathcal{M}_k(\Z_n)$, $M(A)$ is a free
module if and only if $A \in GS(k)$.
\item[(d)] For each $r \in \Z_n^+$, $M(r)$ is a stably free non-free
module if and only if $(r \oplus I_k) \in GS(k+1)$ for some $k$ but $r
\notin GS(1)$.
\end{itemize} \end{prop} 
Note that we may apply the Swan Jakobinski theorem in order to reduce
$(d)$ to the case where $k=1$.

For each unit $r \in \Z_n$, the {\sl Swan module} $(N,r)$ is the
submodule of $\zg$ generated by the elements $N$ and $\hat r$, where
$\hat r$ is any element such that $\varepsilon(\hat r) = r$ mod
$n$. Note that $(N,r) \cong M(r)$. Swan modules, as originally defined
in \cite{Swanper}, are projective and form a well studied subset of
the projective class group of finite groups. We shall only use the
fact that they represent the class of modules $M(r)$ for $r \in \Z_n$
which are projective.

Now let $G$ denote the quaternion group of order $4n$ with $n \geq 6$,
or in general any finite group of period four such that there are
stably free modules $\zg$ which are not free.  We say that weak
cancellation holds for $G$ if all stably free Swan modules are
free. Johnson has shown \cite{Q}, that if weak cancellation holds for
$G$, then one may construct modules which are stably equivalent to
$\pi_2( X_{\mathcal{G}})$ for some presentation $\mathcal{G}$, but
which are {\sl not} isomorphic to $\pi_2( X_{\mathcal{G}}) \oplus
\zg^m$. We refer the interested reader to the recently published Beyl
and Waller \cite{beyl} for an explicit construction of such a module.

For $G$ a 2--group, or for $G$ of order $4p$ with $p$ an odd prime,
Swan has shown that weak cancellation holds \cite[Theorem IV]{bin}.

\section{Appendix}
\begin{prop} For $G=C_n$, the integral group ring $\zg$ has finitely many non-projective maximal ideals.
\end{prop}
\begin{proof} Let $J$ be a maximal ideal of $\zg$. Then $\zg / J = \mathbb{F}$ is a field, which is necessarily finite with ground ring $\Z_p$ for some prime $p$. Define $J \co \zg$ to be the set of elements $r \in \Z$ such that $r \zg \subseteq J$. Clearly, $J \co \zg = (p)$ and so by proposition 7.1 of \cite{Swanproj} $J$ is projective unless $p$ divides the order of $G$. 

Thus, if $J$ is not projective, we may assume that $J \cap \Z = (p)$
for some $p$ dividing $n$, where we consider $\Z$ to be a subring of
$\zg$. Again, since $\zg / J =\mathbb{F}$ is a field, the generator
$x$ for $C_n$ has some minimal polynomial $\omega(x)$ over
$\mathbb{F}$. Then $J$ is necessarily the ideal generated by $p$ and
$\omega(x)$. Since the degree of $\omega(x)$ is less than or equal to
$n$, there are finitely many $\omega(x)$ such that the ideals
$(p,\omega(x))$ are distinct. This completes the proof.
\end{proof}

\begin{prop} The D(2) property holds for $C_{\infty} \times C_{\infty}$. \end{prop}
\begin{proof} If $\Gamma$ is a free abelian group, then every projective $\zG$ module is free, see for example proposition 4.12 of \cite{Lam}. In particular all stably free modules are free.
 We may take the the presentation $\mathcal{G} = \langle x,t | xt=tx
\rangle$ for $C_{\infty} \times C_{\infty} = \Gamma$, leading to the
algebraic complex: \equas \xymatrix{ 0 \ar[r] & \zG
\ar[r]^{\partial_2} & \zG^2 \ar[r]^{\partial_1} & \zG
\ar[r]^{\varepsilon} & 0} \\ \partial_2 = \bem 1-t \\ x-1 \enm \quad 
\quad \partial_1 = (x-1 \quad t-1). \equaf Immediately we deduce that
$\Ext^3(\Z,M)=0$ for any $\zG$--module $M$, and if there exists an
exact sequence: $$0 \lr M \lr F_2 \stackrel{\partial_2}{\lr} F_1 {\lr}
F_0 {\lr} \Z \lr 0$$ with each $F_i$ free, then $M$ is free. This may
be seen to complete the proof by \fullref{red1}. \end{proof} Other
than $S^2$ and $\mathbb{R}P^2$, which are dealt with in \cite{SMD2},
the second homotopy module of any 2--manifold is necessarily zero, and
hence if $\Gamma$ is the fundamental group of {\sl any} surface and
all stably free $\zG$ modules are free, a similar proof shows that the
D(2) property holds for $\Gamma$.

\medskip 
{\bf Acknowledgements}\qua The author wishes to express his gratitude
to his PhD supervisor, Professor F\,E\,A Johnson, who suggested the
problem and who provided valuable insights regarding its eventual
solution. The author also wishes to thank EPSRC for providing
sponsorship.

\bibliographystyle{gtart} \bibliography{link}

\begin{thebibliography}{}
\providecommand\bibmarginpar{\leavevmode\marginpar}
\def\urlstyle#1{{\tt #1}}

\bibitem{Bass}
\textbf{H Bass}, \textbf{M\,P Murthy},
  \href{http://links.jstor.org/sici?sici=0003-486X(196707)2:86:1%3C16:GGAPGO%3%
E2.0.CO%3B2--H} {\emph{Grothendieck groups and {P}icard groups of abelian group
  rings}}, Ann. of Math. $(2)$ 86 (1967) 16--73 \xox{MR}{0219592}

\bibitem{beyl}
\textbf{F\,R Beyl}, \textbf{N Waller},
  \href{http://dx.doi.org/10.2140/agt.2005.5.899} {\emph{A stably free nonfree
  module and its relevance for homotopy classification, case $Q_{28}$}},
  Algebr. Geom. Topol. 5 (2005) 899--910 \xox{MR}{2171797}

\bibitem{browning}
\textbf{W Browning}, \emph{Homotopy types of certain finite C.W. complexes with
  finite fundamental group}, PhD thesis, Cornell University (1978)

\bibitem{CEHom}
\textbf{H Cartan}, \textbf{S Eilenberg}, \emph{Homological algebra}, Princeton
  University Press, Princeton, N. J. (1956) \xox{MR}{0077480}

\bibitem{Fox}
\textbf{R\,H Fox},
  \href{http://links.jstor.org/sici?sici=0003-486X(196005)2:71:3%3C408:FDCVTA%%
3E2.0.CO%3B2--A} {\emph{Free differential calculus V: {T}he {A}lexander
  matrices re-examined}}, Ann. of Math. $(2)$ 71 (1960) 408--422
  \xox{MR}{0111781}

\bibitem{jodie}
\textbf{J\,J A\,M Humphreys}, \emph{Algebraic properties of semi-simple
  lattices and related groups}, PhD thesis, University of London (2005)

\bibitem{Min}
\textbf{F\,E\,A Johnson}, \href{http://dx.doi.org/10.1017/S0305004102006011}
  {\emph{Explicit homotopy equivalences in dimension two}}, Math. Proc.
  Cambridge Philos. Soc. 133 (2002) 411--430 \xox{MR}{1919714}

\bibitem{SMD2}
\textbf{F\,E\,A Johnson}, \emph{Stable modules and the $D(2)$--problem}, London
  Mathematical Society Lecture Note Series 301, Cambridge University Press,
  Cambridge (2003) \xox{MR}{2012779}

\bibitem{Q}
\textbf{F\,E\,A Johnson},
  \href{http://www.ams.org/jourcgi/jour-getitem?pii=S0002993903070680}
  {\emph{Minimal 2--complexes and the $\rm D(2)$--problem}}, Proc. Amer. Math.
  Soc. 132 (2004) 579--586 \xox{MR}{2022384}

\bibitem{Lam}
\textbf{T\,Y Lam}, \emph{Serre's conjecture}, Springer, Berlin (1978)
  \xox{MR}{0485842}

\bibitem{lat}
\textbf{M\,P Latiolais},
  \href{http://links.jstor.org/sici?sici=0002-9947(198602)293:2%3C655:SHTOF2%3%
E2.0.CO%3B2--X} {\emph{Simple homotopy type of finite 2--complexes with finite
  abelian fundamental group}}, Trans. Amer. Math. Soc. 293 (1986) 655--662
  \xox{MR}{816317}

\bibitem{Mac}
\textbf{S Mac~Lane}, \emph{Homology}, Classics in Mathematics, Springer, Berlin
  (1995) \xox{MR}{1344215}

\bibitem{wajid}
\textbf{W Mannan}, \emph{PhD Thesis}, in preparation

\bibitem{Milnor}
\textbf{J Milnor}, \emph{Introduction to algebraic $K$--theory}, Princeton
  University Press, Princeton, N.J. (1971) \xox{MR}{0349811}

\bibitem{Swanproj}
\textbf{R\,G Swan},
  \href{http://links.jstor.org/sici?sici=0003-486X(196005)2:71:3%3C552:IRAPM%3%
E2.0.CO%3B2--F} {\emph{Induced representations and projective modules}}, Ann.
  of Math. $(2)$ 71 (1960) 552--578 \xox{MR}{0138688}

\bibitem{Swanper}
\textbf{R\,G Swan},
  \href{http://links.jstor.org/sici?sici=0003-486X(196009)2:72:2%3C267:PRFFG%3%
E2.0.CO%3B2--X} {\emph{Periodic resolutions for finite groups}}, Ann. of Math.
  $(2)$ 72 (1960) 267--291 \xox{MR}{0124895}

\bibitem{bin}
\textbf{R\,G Swan}, \emph{Projective modules over binary polyhedral groups}, J.
  Reine Angew. Math. 342 (1983) 66--172 \xox{MR}{703486}

\bibitem{fc}
\textbf{C\,T\,C Wall},
  \href{http://links.jstor.org/sici?sici=0003-486X(196501)2:81:1%3C56:FCFC%3E2%
.0.CO%3B2--F} {\emph{Finiteness conditions for $\mathrm{CW}$--complexes}}, Ann.
  of Math. $(2)$ 81 (1965) 56--69 \xox{MR}{0171284}

\end{thebibliography}

\end{document}